\documentclass[12pt, twoside, leqno]{article}

\usepackage{amsmath,amsthm}
\usepackage{amssymb}
\usepackage{amsfonts,bm,float,inputenc, mathrsfs, bbm, cite, hyperref}
\usepackage{enumitem}

\usepackage{graphicx}

\usepackage[T1]{fontenc}

\newtheorem{theorem}{Theorem}[section]

\newtheorem{lemma}[theorem]{Lemma}
\newtheorem{proposition}[theorem]{Proposition}

\newtheorem{conjecture}[theorem]{Conjecture}
\newtheorem{question}[theorem]{Question}

\theoremstyle{definition}
\newtheorem{definition}[theorem]{Definition}
\newtheorem{remark}[theorem]{Remark}

\numberwithin{equation}{section}

\begin{document}


\baselineskip=17pt


\title{Pointwise Ergodic Theorems for Higher Levels of Mixing}

\author{Sohail Farhangi\\
Department of Mathematics\\ 
The Ohio State University\\
Mathematics Building 430\\
Columbus, Ohio 43220\\
E-mail: sohail.farhangi@gmail.com}

\date{}

\maketitle


\renewcommand{\thefootnote}{}



\renewcommand{\thefootnote}{\arabic{footnote}}
\setcounter{footnote}{0}


\begin{abstract}
We prove strengthenings of the Birkhoff Ergodic Theorem for weakly mixing and strongly mixing measure preserving systems. We show that our pointwise theorem for weakly mixing systems is strictly stronger than the Wiener-Wintner Theorem. We also show that our pointwise Theorems for weakly mixing and strongly mixing systems characterize weakly mixing systems and strongly mixing systems respectively.
\end{abstract}

\section{Introduction}

\hskip 4mm In this section we establish the notation that we use, review some known results related to the Birkhoff Pointwise Ergodic Theorem and state the main theorems of this paper.

Whenever we discuss a measure preserving system (m.p.s.) $(X,\mathscr{B},\mu,T)$, $X$ will be a measurable space, $\mathscr{B}$ will be a $\sigma$-algebra of subsets of $X$, $\mu$ will be a probability measure on $(X,\mathscr{B})$ and $T:X\rightarrow X$ will be a measurable transformation satisfying $\mu(A) = \mu(T^{-1}A)$ for all $A \in \mathscr{B}$. When we work with the Hilbert space $L^2(X,\mu)$, we will let $U:L^2(X,\mu)\rightarrow L^2(X,\mu)$ denote the unitary operator given by $U(f) = f\circ T$. When we work with a m.p.s. of the form $([0,1],\mathscr{B},\mu,T)$, we will assume that $\mathscr{B}$ is the completion of the Borel $\sigma$-algebra. We say that two sequences of complex numbers $(x_n)_{n = 1}^{\infty}$ and $(y_n)_{n = 1}^{\infty}$ are orthogonal if

\begin{equation}
    \lim_{N\rightarrow\infty}\frac{1}{N}\sum_{n = 1}^Nx_n\overline{y_n} = 0.
\end{equation}

Now let us recall the levels of the ergodic hierarchy of mixing that will be used throughout this paper.

\begin{definition}
\label{EHOfMixing}
Let $(X,\mathscr{B},\mu,T)$ be a m.p.s.

\begin{itemize}
\item $(X,\mathscr{B},\mu,T)$ is {\bf ergodic} if for every $A \in \mathscr{B}$ satisfying $\mu(A) = \mu(T^{-1}A)$ we have $\mu(A) \in \{0,1\}$.

\item $(X,\mathscr{B},\mu,T)$ is {\bf weakly mixing} if for every $A,B \in \mathscr{B}$ we have

\begin{equation}
\lim_{N\rightarrow\infty}\frac{1}{N}\sum_{n = 1}^N|\mu(A\cap T^{-n}B)-\mu(A)\mu(B)| = 0.
\end{equation}

\item $(X,\mathscr{B},\mu,T)$ is {\bf strongly mixing} if for every $A,B \in \mathscr{B}$ we have

\begin{equation}
\lim_{n\rightarrow\infty}\mu(A\cap T^{-n}B) = \mu(A)\mu(B).
\end{equation}
\end{itemize}
\end{definition}

\begin{theorem}[Birkhoff] Let $(X,\mathscr{B},\mu,T)$ be a  m.p.s., and let $f \in L^1([0,1],\mu)$. For a.e. $x \in X$, we have
\label{BET}
\begin{equation}
\lim_{N\rightarrow\infty}\frac{1}{N}\sum_{n = 1}^Nf(T^nx) = f^*(x),
\end{equation}

\noindent where $f^*(x) \in L^1(X,\mu)$ is such that $f^*(Tx) = f^*(x)$ for a.e. $x \in X$ and $\int_Af^*d\mu = \int_Afd\mu$ for every $A \in \mathscr{B}$ satisfying $A = T^{-1}(A)$. In particular, if $T$ is ergodic, then for a.e. $x \in X$ we have

\begin{equation}
\lim_{N\rightarrow\infty}\frac{1}{N}\sum_{n = 1}^Nf(T^nx) = \int_Xfd\mu. 
\end{equation}
\end{theorem}

The Birkhoff Pointwise Ergodic Theorem can be interpreted as follows. Given an ergodic m.p.s. $(X,\mathscr{B},\mu,T)$ and $f \in L^1(X,\mu)$ satisfying $\int_Xfd\mu = 0$, the sequence $\left(f(T^nx)\right)_{n = 1}^{\infty}$ is orthogonal to the constant sequence $(1)_{n = 1}^{\infty}$. The Wiener-Wintner Theorem is a generalization of the Birkhoff Pointwise Ergodic Theorem for weakly mixing systems and has a similar interpretation.

\begin{theorem}[Wiener-Wintner]\label{WWT} Let $(X,\mathscr{B},\mu,T)$ be a m.p.s. and let $f \in L^1(X,\mu)$. There exists $X' \in \mathscr{B}$ with $\mu(X') = 1$, such that for every $x \in X'$ and any $\lambda \in \mathbb{C}$ with $|\lambda| = 1$ the limit

\begin{equation}
\lim_{N\rightarrow\infty}\frac{1}{N}\sum_{n = 1}^Nf(T^nx)\lambda^n
\end{equation}

\noindent exists. Furthermore, if $T$ is weakly mixing, then

\begin{equation}
\lim_{N\rightarrow\infty}\frac{1}{N}\sum_{n = 1}^Nf(T^nx)\lambda^n = \begin{cases}
                                                                        0 & \text{if }\lambda \neq 1\\
																		\int_Xfd\mu & \text{if } \lambda = 1
			     													 \end{cases}.
\end{equation}
\end{theorem}


\begin{definition} A bounded sequence of complex numbers $(x_n)_{n = 1}^{\infty}$ is \noindent{\bf Besicovitch Almost Periodic} if for each $\epsilon > 0$ there exists a trigonometric polynomial $P_{\epsilon}$ for which

\begin{equation}
\overline{\lim_{N\rightarrow\infty}}\frac{1}{N}\sum_{n = 1}^N|x_n-P_{\epsilon}(n)| < \epsilon.
\end{equation}
\end{definition}


The Wiener-Wintner Theorem can be interpreted as follows. Given a weakly mixing m.p.s. $(X,\mathscr{B},\mu,T)$ and $f \in L^1(X,\mu)$ satisfying $\int_Xfd\mu = 0$, the sequence $\left(f(T^nx)\right)_{n = 1}^{\infty}$ is orthogonal to any Besicovitch Almost Periodic Sequence $(y_n)_{n = 1}^{\infty}$. Theorem \ref{FirstMainTheorem} is one of the main results of this paper and is a generalization of the Wiener-Wintner Theorem when interpreted in a similar fashion. We require some more preliminaries before we can state Theorem \ref{FirstMainTheorem}.


\begin{definition} Let $(x_n)_{n = 1}^{\infty} \subseteq \mathbb{C}$ satisfy

\begin{equation}
\overline{\lim_{N\rightarrow\infty}}\frac{1}{N}\sum_{n = 1}^N|x_n| < \infty.
\end{equation}

\noindent $(x_n)_{n = 1}^{\infty}$ is a \noindent{\bf weakly mixing sequence} if it satisfies the following condition.

Suppose that for a bounded sequence $(y_n)_{n = 1}^{\infty} \subseteq \mathbb{C}$ there is a strictly increasing sequence $(N_q)_{q = 1}^{\infty} \subseteq \mathbb{N}$ for which

\begin{equation}
\ell(h) := \lim_{q\rightarrow\infty}\frac{1}{N_q}\sum_{n = 1}^{N_q}x_{n+h}\overline{y_n}\text{ exists for every }h \in \mathbb{N}.
\end{equation}
Then we have
\begin{equation}
\lim_{H\rightarrow\infty}\frac{1}{H}\sum_{h = 1}^H|\ell(h)| = 0.
\end{equation}
\end{definition}


\begin{theorem}[cf. Theorem \ref{ActualyFirstMainTheorem} in Section 2] Let $(X,\mathscr{B},\mu,T)$ be a weakly mixing m.p.s. and let $f \in L^1(X,\mu)$ satisfy $\int_Xfd\mu = 0$. For a.e. $x \in X$, $\left(f(T^nx)\right)_{n = 1}^{\infty}$ is a weakly mixing sequence.
\label{FirstMainTheorem}
\end{theorem}


\begin{definition}[cf. Definition 3.13 in \cite{TheErdosSumsetPaper}] A bounded sequence of complex numbers $(x_n)_{n = 1}^{\infty}$ is \noindent{\bf compact} if for any $\epsilon > 0$, there exists $K \in \mathbb{N}$ such that

\begin{equation}
\underset{m \in \mathbb{N}}{\text{sup }}\underset{1 \le k \le K}{\text{min }}\overline{\lim_{N\rightarrow\infty}}\frac{1}{N}\sum_{n = 1}^N|x_{n+m}-x_{n+k}|^2 < \epsilon.
\end{equation}
\end{definition}


\begin{lemma}[cf. Lemma 3.19 of \cite{TheErdosSumsetPaper}] If $(w_n)_{n = 1}^{\infty} \subseteq \mathbb{C}$ is a weakly mixing sequence and $(c_n)_{n = 1}^{\infty}$ is a compact sequence, then

\begin{equation}
\lim_{N\rightarrow\infty}\frac{1}{N}\sum_{n = 1}^Nw_n\overline{c_n} = 0.
\end{equation}
\label{WeakPerpCom}
\end{lemma}


Lemma \ref{WeakPerpCom} is what allows us to interpret Theorem \ref{FirstMainTheorem} in a similar fashion to the Birkhoff Pointwise Ergodic Theorem and the Wiener-Wintner Theorem. Theorem \ref{FirstMainTheorem} asserts that for a weakly mixing m.p.s. $(X,\mathscr{B},\mu,T)$ and $f \in L^1(X,\mu)$ satisfying $\int_Xfd\mu = 0$, the sequence $\left(f(T^nx)\right)_{n = 1}^{\infty}$ is orthogonal to any compact sequence  $(y_n)_{n = 1}^{\infty}$. We will see in section $2$ that the class of compact sequences is strictly larger than the class of Besicovitch Almost Periodic Sequences and this is why Theorem \ref{FirstMainTheorem} is more general than the Wiener-Wintner Theorem. The next main result of this paper is an analogue of Theorem \ref{FirstMainTheorem} for strongly mixing measure preserving systems.

\begin{definition}
Let $(x_n)_{n = 1}^{\infty} \subseteq \mathbb{C}$ satisfy

\begin{equation}
\overline{\lim_{N\rightarrow\infty}}\frac{1}{N}\sum_{n = 1}^N|x_n| < \infty.
\end{equation}

\noindent $(x_n)_{n = 1}^{\infty}$ is a \noindent{\bf strongly mixing sequence} if it satisfies the following condition.

Suppose that for a bounded sequence $(y_n)_{n = 1}^{\infty} \subseteq \mathbb{C}$ there is a strictly increasing sequence $(N_q)_{q = 1}^{\infty} \subseteq \mathbb{N}$ for which

\begin{equation}
\ell(h) := \lim_{q\rightarrow\infty}\frac{1}{N_q}\sum_{n = 1}^{N_q}x_{n+h}\overline{y_n}\text{ exists for every }h \in \mathbb{N}.
\end{equation}
Then we have
\begin{equation}
\lim_{h\rightarrow\infty}|\ell(h)| = 0.
\end{equation}
\end{definition}


\begin{theorem}[cf. Theorem \ref{mainStrongMixingTheorem} of Section 3] Let $(X,\mathscr{B},\mu,T)$ be a strongly mixing m.p.s. and let $f \in L^1(X,\mu)$ satisfy $\int_Xfd\mu = 0$. For a.e. $x \in X$, $\left(f(T^nx)\right)_{n = 1}^{\infty}$ is a strongly mixing sequence.
\label{SecondMainTheorem}
\end{theorem}


Proposition \ref{TooTrivial} gives a partial converse to the Birkhoff Pointwise Ergodic Theorem. Proposition \ref{TooTrivial} is well known and an easy consequence of the Dominated Convergence Theorem but we state it anyways so that we can put some results of this paper in context.


\begin{proposition} Let $([0,1],\mathscr{B},\mu,T)$ be a m.p.s. If for every $f \in L^{\infty}([0,1],\mu)$, there exists $A_f \in \mathscr{B}$ such that $\mu(A_f) = 1$ and for every $x \in A_f$ we have

\begin{equation}
\label{BirkConverse}
\lim_{N\rightarrow\infty}\frac{1}{N}\sum_{n = 1}^Nf(T^nx) = \int_Xfd\mu,
\end{equation}

\noindent then $T$ is ergodic.
\label{TooTrivial}
\end{proposition}


Proposition \ref{Easy}  is a converse to Theorem \ref{FirstMainTheorem} in the same fashion that Proposition \ref{TooTrivial} is a converse to the Birkhoff Ergodic Theorem.


\begin{proposition} Let $(X,\mathscr{B},\mu,T)$ be a m.p.s. If for every $f \in L^{\infty}(X,\mu)$ with $\int_Xfd\mu = 0$ there exist a set $A_f \subseteq X$ satisfying $\mu(A_f) = 1$ and for every $x \in A_f$ we have that $\left(f(T^nx)\right)_{n = 1}^{\infty}$ is a weakly mixing sequence, then $T$ is weakly mixing.
\label{Easy}
\end{proposition}


Proposition \ref{Easy} should not come as a surprise since the Wiener-Wintner Theorem is already strong enough to characterize weakly mixing systems. To see this let us recall that an ergodic m.p.s. $(X,\mathscr{B},\mu,T)$ is weakly mixing if and only if $L^2(X,\mu)$ has no non-constant eigenfunctions with respect to $U$. Now let $(X,\mathscr{B},\mu,T)$ be an ergodic m.p.s. that is not weakly mixing and let $f \in L^2(X,\mu)$ be a non-constant eigen function corresponding to the eigenvalue $\lambda$. We see that

\begin{multline}
\lim_{H\rightarrow\infty}\frac{1}{H}\sum_{h = 1}^H|\lim_{N\rightarrow\infty}\frac{1}{N}\sum_{n = 1}^Nf(T^{n+h}x)\lambda^{-n}| \\ = \lim_{H\rightarrow\infty}\frac{1}{H}\sum_{h = 1}^H|f(T^hx)| = \int_X|f|d\mu \neq 0
\end{multline}

\noindent for a.e. $x \in X$. Since any m.p.s. $(X,\mathscr{B},\mu,T)$ satisfying the hypothesis of Proposition \ref{Easy} also satisfies the hypothesis of Proposition \ref{TooTrivial} we see that $(X,\mathscr{B},\mu,T)$ is ergodic, so $(X,\mathscr{B},\mu,T)$ is in fact weakly mixing.

Theorem \ref{StrongMixingConverse} is a converse to Theorem \ref{SecondMainTheorem} in the same way that Proposition \ref{Easy} is a converse for Theorem \ref{FirstMainTheorem}. The proof of Theorem \ref{StrongMixingConverse} can easily be adjusted to give a direct proof of Proposition \ref{Easy} as well, but we will delay the proof of Theorem \ref{StrongMixingConverse} until section 3.


\begin{theorem}[cf. Theorem \ref{ActualStrongMixingConverse} in Section 3] Let $([0,1],\mathscr{B},\mu,T)$ be a m.p.s. If for every $f \in L^{\infty}([0,1],\mu)$ with $\int_Xfd\mu = 0$ there exist a set $A_f \subseteq X$ satisfying $\mu(A_f) = 1$ and for every $x \in A_f$ we have that $\left(f(T^nx)\right)_{n = 1}^{\infty}$ is a strongly mixing sequence, then $T$ is strongly mixing.
\label{StrongMixingConverse}
\end{theorem}
\section{A Pointwise Theorem for Weakly Mixing Systems}

\hskip 4mm In this section we prove Theorem \ref{ActualyFirstMainTheorem}, which is Theorem \ref{FirstMainTheorem} for standard measure preserving systems. In section $4$ we will show how to deduce Theorem \ref{FirstMainTheorem} from \ref{ActualyFirstMainTheorem}. We will also discuss some potential generalizations of Theorem \ref{FirstMainTheorem}.

Before we begin proving Theorem \ref{ActualyFirstMainTheorem} we will show that the class of compact sequences is strictly larger than that of Besicovitch Almost Periodic Sequences. In Example 3.21 of \cite{TheErdosSumsetPaper} a compact sequence orthogonal to any Besicovitch Almost Periodic Sequence was constructed. In Lemma \ref{Example} we give a different class of compact sequences that are orthogonal to all Besicovitch Almost Periodic Sequences.


\begin{lemma} Let $\alpha \in \mathbb{T}$ be irrational, and let $(x_n)_{n = 1}^{\infty}$ be the sequence 

$$\alpha, 2\alpha, 2\alpha, 3\alpha, 3\alpha, 3\alpha, \cdots,$$

\noindent which can also be expressed as $x_n = m\alpha$ for $\binom{m}{2} < n \le \binom{m+1}{2}$. If $f \in C(\mathbb{T})$ satisfies $\int_{\mathbb{T}}fdm = 0$, then $\left(f(x_n)\right)_{n = 1}^{\infty}$ is orthogonal to any Besicovitch Almost Periodic Sequence.
\label{Example}
\end{lemma}


\noindent{\it Proof.} Using standard approximation arguments, we see that it suffices to show that for any $k \in \mathbb{N}$ the sequence $\left(e^{2\pi i kx_n}\right)_{n = 1}^{\infty}$ is orthogonal to any Besicovitch Almost Periodic Sequence. We note that for any $\beta \in \mathbb{T}\setminus\{0\}$ and $k \in \mathbb{Z}\setminus\{0\}$, we have

\begin{equation}
\overline{\lim_{N\rightarrow\infty}}|\frac{1}{N}\sum_{n = 1}^Ne^{2\pi ikx_n}e^{-2\pi in\beta}|
\end{equation}

\begin{equation}
= \overline{\lim_{N\rightarrow\infty}}|\binom{N}{2}^{-1}\sum_{n = 1}^Ne^{-2\pi i\binom{n}{2}\beta}\sum_{m = 1}^ne^{2\pi ikn\alpha}e^{-2\pi im\beta}|
\end{equation}

\begin{equation}
= \overline{\lim_{N\rightarrow\infty}}|\binom{N}{2}^{-1}\sum_{n = 1}^Ne^{2\pi i(n\alpha-\binom{n}{2}\beta)}\frac{e^{-2\pi i(n+1)\beta}-e^{-2\pi i\beta}}{1-e^{-2\pi i\beta}}| 
\end{equation}

\begin{equation}
\le \overline{\lim_{N\rightarrow\infty}}\binom{N}{2}^{-1}\frac{N}{|1-e^{-2\pi i\beta}|} = 0.
\end{equation}

\noindent Similarly, we see that 

\begin{equation}
\overline{\lim_{N\rightarrow\infty}}|\frac{1}{N}\sum_{n = 1}^Ne^{2\pi ikx_n}\overline{1}| = \overline{\lim_{N\rightarrow\infty}}|\binom{N}{2}^{-1}\sum_{n = 1}^Nne^{2\pi ikn\alpha}|
\end{equation}

\begin{equation}
= \overline{\lim_{N\rightarrow\infty}}|\frac{1}{n}\sum_{n = 1}^Ne^{2\pi ib\alpha}| = 0.
\end{equation}

We now see that a standard approximation argument shows that $(e^{2\pi ikx_n})_{n = 1}^{\infty}$ is orthogonal to any Besicovitch Almost Periodic Sequence.

To see that $\left(f(x_n)\right)_{n = 1}^{\infty}$ is a compact sequence it suffices to note that

\begin{equation}
   \underset{m \in \mathbb{N}}{\text{sup }}\underset{1 \le k \le K}{\text{min }}\overline{\lim_{N\rightarrow\infty}}\frac{1}{N}\sum_{n = 1}^N|x_{n+m}-x_{n+k}|^2
\end{equation}
\begin{equation}
\le \underset{m \ge K+1}{\text{sup }}\underset{1 \le k \le K}{\text{min }}\sum_{j = k}^{m-1} \overline{\lim_{N\rightarrow\infty}}\frac{1}{N}\sum_{n = 1}^N|x_{n+j+1}-x_{n+j}|^2 
\end{equation}

\begin{equation}
\pushQED{\qed}
    = \underset{m \ge k+1}{\text{sup }}\underset{1 \le k \le K}{\text{min }}\sum_{j = k}^{m-1}0 = 0. \qedhere
\popQED
\end{equation}


We see that in the proof of Lemma \ref{Example}, the sequence $(n\alpha)_{n = 1}^{\infty}$ could be replaced by any uniformly distributed sequence in $[0,1]$, and $(\binom{m}{2})_{m = 1}^{\infty}$ could be replaced by any subpolynomial sequence of integers $(k_m)_{m = 1}^{\infty}$ satisfying $\lim_{m\rightarrow\infty}k_{m+1}-k_m = \infty$, so Lemma \ref{Example} can be used to construct many examples of compact sequences that are not Besicovitch Almost Periodic. 

Lemma \ref{Representation} is what will allow us to view bounded sequences of complex numbers as an element of some Hilbert space. 

\begin{lemma}[cf. Lemma 3.26 of \cite{TheErdosSumsetPaper}] Let $(a_n)_{n = 1}^{\infty} \subseteq \mathbb{C}$ be bounded. Then there exists a compact metric space $X$, a continuous map $S:X\rightarrow X$, a function $F \in C([0,1])$, and a point $x \in X$ with a dense orbit under $S$ such that $a_n = F(S^n(x))$ for all $n \in \mathbb{N}$. 
\label{Representation}
\end{lemma}

Lemma \ref{keylemma} allows us to convert correlations of sequences to inner products of functions in a Hilbert space, which will then permit us to use Hilbert space Theory to analyze our correlations. Before stating Lemma \ref{keylemma}, let us recall the definition of a generic point.

\begin{definition} Given an ergodic m.p.s. $([0,1],\mathscr{B},\mu,T)$, $x \in [0,1]$ is {\bf generic} if for every $f \in C([0,1])$ we have

\begin{equation}
\lim_{N\rightarrow\infty}\frac{1}{N}\sum_{n = 1}^Nf(T^nx) = \int_{[0,1]}fd\mu.
\end{equation}
\end{definition}

\begin{lemma}\label{keylemma} Let $([0,1],\mathscr{B},\mu,T)$ be an ergodic m.p.s. Let $(y_n)_{n = 1}^{\infty}$ be a bounded sequence of complex numbers. Let $f \in C([0,1])$ be arbitrary and let $x \in X$ be a generic point. Let $U:L^2([0,1],\mu)\rightarrow L^2([0,1],\mu)$ be the unitary operator induced by $T$. Let $(N_q)_{q = 1}^{\infty} \subseteq \mathbb{N}$ be any sequence for which

\begin{equation}
\lim_{q\rightarrow\infty}\frac{1}{N_q}\sum_{n = 1}^{N_q}f(T^{n+h}x)\overline{y_n}
\end{equation}

\noindent exists for each $h \in \mathbb{N}$. Then there exists $g \in L^2([0,1],\mu)$ for which

\begin{equation}
\lim_{q\rightarrow\infty}\frac{1}{N_q}\sum_{n = 1}^{N_q}f(T^{n+h}x)\overline{y_n} = \langle U^hf,g\rangle_{\mu}.
\end{equation}
\end{lemma}


\noindent{\it Proof.} By Lemma \ref{Representation}, let $(Y,d)$ be a compact metric space, let $\xi \in C(Y)$, $y \in Y$, and $S:Y\rightarrow Y$ be such that $y_n = \xi(S^ny)$. Let $\nu$ be any weak$^*$ limit point of the sequence 

\begin{equation}
\left(\frac{1}{N_q}\sum_{k = 1}^{N_q}\delta_{T^nx,S^ny}\right)_{q = 1}^{\infty},
\end{equation}
and let $(M_q)_{q = 1}^{\infty} \subseteq \mathbb{N}$ be such that

\begin{equation}
\nu = \lim_{q\rightarrow\infty}\frac{1}{M_q}\sum_{n = 1}^{M_q}\delta_{T^nx,S^ny}.
\end{equation}

Let $\tilde{f},\tilde{\xi} \in L^2([0,1]\times Y,\nu)$ be given by $\tilde{f}(x,y) = f(x)$ and $\tilde{\xi}(x,y) = \xi(y)$. Let $V:L^2([0,1]\times Y,\nu)\rightarrow L^2([0,1]\times Y,\nu)$ be the unitary operator induced by $T\times S$. We note that if $h \in L^2([0,1]\times Y,\nu)$ is such that $h(x,y) = k(x)$ for some $k \in L^2([0,1],\mu)$, then the genericity of $x$ gives us

\begin{equation}
\int_{[0,1]\times Y}hd\nu = \int_{[0,1]}kd\mu.
\end{equation}

Let $\tilde{\mu}$ be the probability measure on $([0,1]\times Y,\mathscr{B}\otimes\mathscr{A})$ given by $\tilde{\mu}(B\times A) = \mu(B)\mathbbm{1}_A(y)$ for all $B \in \mathscr{B}$ and $A \in \mathscr{A}$. Since we may identify $L^2([0,1]\times Y,\tilde{\mu})$ with the functions in $L^2([0,1]\times Y,\nu)$ of the form $h(x,y) = k(x)$, let $P:L^2([0,1]\times Y,\nu)\rightarrow L^2([0,1]\times Y,\tilde{\mu})$ denote the orthogonal projection. Let $\tilde{g} = P\tilde{\xi}$, and let $g \in L^2([0,1],\mu)$ be such that $\tilde{g}(x,y) = g(x)$. We now see that for any $h \in \mathbb{N}$ we have

\begin{equation}
\lim_{q\rightarrow\infty}|\frac{1}{N_q}\sum_{n = 1}^{N_q}f(T^{n+h}x)\overline{y_n}| = \lim_{q\rightarrow\infty}|\frac{1}{M_q}\sum_{n = 1}^{M_q}\tilde{f}(T^{n+h}x)\overline{\tilde{\xi}(S^ny)}|
\end{equation}

\begin{equation}
\pushQED{\qed}
= |\int_{[0,1]\times Y} V^hf\overline{\tilde{\xi}}d\nu| = |\langle V^h\tilde{f}, \tilde{\xi}\rangle_{\nu}| = |\langle V^h\tilde{f},\tilde{g}\rangle_{\nu}| = |\langle U^hf,g\rangle_{\mu}|. \qedhere
\popQED
\end{equation}


\begin{lemma} Let $([0,1],\mathscr{B},\mu,T)$ be a weakly mixing m.p.s. For each generic point $x \in [0,1]$ and each $f \in C([0,1])$ with $\int_{[0,1]}fd\mu = 0$, $\left(f(T^nx)\right)_{n = 1}^{\infty}$ is a weakly mixing sequence.
\label{WeakMixingForC(X)}
\end{lemma}


\noindent{\it Proof.} Let $f \in C([0,1])$ be arbitrary and let $(y_n)_{n = 1}^{\infty} \subseteq \mathbb{C}$ be bounded. Let $(N_q)_{q = 1}^{\infty} \subseteq \mathbb{N}$ be any sequence for which

\begin{equation}
\lim_{q\rightarrow\infty}\frac{1}{N_q}\sum_{n = 1}^{N_q}f(T^{n+h}x)\overline{y_n}
\end{equation}

\noindent exists for every $h \in \mathbb{N}$. By Lemma \ref{keylemma}, let $g \in L^2([0,1],\mu)$ be such that

\begin{equation}
\lim_{q\rightarrow\infty}\frac{1}{N_q}\sum_{n = 1}^{N_q}f(T^{n+h}x)\overline{y_n} = \langle U^hf,g\rangle.
\end{equation}

Letting $U:L^2([0,1],\mu)\rightarrow L^2([0,1],\mu)$ denote the unitary operator induced by $T$, we see that $U$ is weakly mixing. It follows that

\begin{equation}
\lim_{H\rightarrow\infty}\frac{1}{H}\sum_{h = 1}^H|\lim_{q\rightarrow\infty}\frac{1}{N_q}\sum_{n = 1}^{N_q}f(T^{n+h}x)\overline{y_n}|
\end{equation}

\begin{equation}
\pushQED{\qed}
= \lim_{H\rightarrow\infty}\frac{1}{H}\sum_{h = 1}^H|\langle U^hf,g\rangle| = 0. \qedhere
\popQED
\end{equation}


\begin{theorem} Let $([0,1],\mathscr{B},\mu,T)$ be a weakly mixing m.p.s. and let $f \in L^1([0,1],\mu)$ satisfy $\int_{[0,1]}fd\mu = 0$. For a.e. $x \in [0,1]$, $\left(f(T^nx)\right)_{n = 1}^{\infty}$ is a weakly mixing sequence.
\label{ActualyFirstMainTheorem}
\end{theorem}


\noindent{\it Proof.} Let $\epsilon > 0$ be arbitrary, and let $g \in C([0,1])$ be such that $||f-g||_1 \le \epsilon$ and $\int_{[0,1]}gd\mu = 0$. Let $X \subseteq [0,1]$ be a set of full measure for which the Birkhoff Pointwise Ergodic Theorem holds for $|f-g|$ along every $x \in X$. We now see that

\begin{equation}
\lim_{N\rightarrow\infty}\frac{1}{N}\sum_{n = 1}^N|f(T^nx)-g(T^nx)| = \int_{[0,1]}|f-g|d\mu = ||f-g||_1 < \epsilon.
\end{equation}

Now let $(y_n)_{n = 1}^{\infty} \subseteq \mathbb{C}$ be uniformly bounded in norm by $1$. Since a.e. $x \in X$ is a generic point we may use Lemma \ref{WeakMixingForC(X)} to further refine $X$ to another set of full measure $X'$, such that for every $x \in X'$, $\left(g(T^nx)\right)_{n = 1}^{\infty}$ is a weakly mixing sequence. We see that for any $x \in X'$, we have

\begin{multline}
\lim_{H\rightarrow\infty}\frac{1}{H}\sum_{h = 1}^H\overline{\lim_{N\rightarrow\infty}}|\frac{1}{N}\sum_{n = 1}^Nf(T^{n+h}x)\overline{y_n}| \\ \le \lim_{H\rightarrow\infty}\frac{1}{H}\sum_{h = 1}^H\overline{\lim_{N\rightarrow\infty}}|\frac{1}{N}\sum_{n = 1}^Ng(T^{n+h}x)\overline{y_n}| +\epsilon = \epsilon.
\label{Modification}
\end{multline}

\noindent Since $\epsilon > 0$ was arbitrary, we see that $\left(f(T^nx)\right)_{n = 1}^{\infty}$ is a weakly mixing sequence. \qed


One way in which we can try to generalize Theorem \ref{ActualyFirstMainTheorem} is motivated by Theorem \ref{BourgainThm} which is due to Bourgain.


\begin{theorem}[Theorem 1 in \cite{Bourgain'sPolynomialTheorem}] Let $(X,\mathscr{B},\mu,T)$ be a m.p.s. and let $p(x)$ be a polynomial with integer coefficients. If $f \in L^r(X,\mu)$ for some $r > 1$, then

\begin{equation}
\lim_{N\rightarrow\infty}\frac{1}{N}\sum_{n = 1}^Nf(T^{p(n)}x)
\end{equation}

\noindent exists for a.e. $x \in X$. Furthermore, if $T$ is weakly mixing, then

\begin{equation}
\lim_{N\rightarrow\infty}\frac{1}{N}\sum_{n = 1}^Nf(T^{p(n)}x) = \int_Xfd\mu
\end{equation}

\noindent for a.e. $x \in X$.
\label{BourgainThm}
\end{theorem}


Theorem \ref{BourgainThm} shows us that the Birkhoff Ergodic Theorem holds along polynomial subsequences if the underlying transformation $T$ is assumed to be weakly mixing. This naturally leads us to ask if Theorem \ref{ActualyFirstMainTheorem} holds for polynomial subsequences.


\begin{question} If $(X,\mathscr{B},\mu,T)$ is a weakly mixing m.p.s., $p(x)$ a polynomial with integer coefficients and $f \in L^r(X,\mu)$ with $r > 1$ is such that $\int_Xfd\mu = 0$, then is $\left(f(T^{p(n)}x)\right)_{n = 1}^{\infty}$ a weakly mixing sequence for a.e. $x \in X$?
\label{PolynomialQuestion}
\end{question}


We see that if $(x_n)_{n = 1}^{\infty}$ is a weakly mixing sequence of complex numbers, and $(y_n)_{n = 1}^{\infty}$ is another sequence of complex numbers for which $d(\{n \in \mathbb{N}\ |\ x_n \neq y_n\})\footnote{For $A \subseteq \mathbb{N}$ the {\bf natural density} of $A$ is given by $$d(A) := \lim_{N\rightarrow\infty}\frac{1}{N}|A\cap[1,N]|$$ provided that the limit exists.} = 0$, then $(y_n)_{n = 1}^{\infty}$ is also weakly mixing. In particular, if $(x_n)_{n = 1}^{\infty}$ is weakly mixing, and $(y_n)_{n = 1}^{\infty}$ is given by $y_n = x_n$ when $n$ is not a square and $y_n = 1$ when $n$ is a square, then $(y_n)_{n = 1}^{\infty}$ is also a weakly mixing sequence, but $(y_{n^2})_{n = 1}^{\infty}$ is the constant sequence, so Question \ref{PolynomialQuestion} does not follow immediately from Theorem \ref{ActualyFirstMainTheorem}. It is well known that if $(X,\mathscr{B},\mu,T)$ is a weakly mixing system, then for any $f,g \in L^2(X,\mu)$ we have

\begin{equation}
\lim_{N\rightarrow\infty}\frac{1}{N}\sum_{n = 1}^N|\langle U^{n^2}f,g\rangle-\int_Xfd\mu\int_Xgd\mu| = 0.
\end{equation}

Combining this fact with Lemma \ref{keylemma} allows us to prove Proposition \ref{PolynomialProp}, but does not immediately help us resolve Question \ref{PolynomialQuestion}.


\begin{proposition} Let $(X,\mathscr{B},\mu,T)$ be a weakly mixing m.p.s. and let $f \in L^1(X,\mu)$ satisfy $\int_{X}fd\mu = 0$. For a.e. $x \in X$, any bounded sequence of complex numbers $(y_n)_{n = 1}^{\infty}$, and any strictly increasing sequence of natural numbers $(N_q)_{q = 1}^{\infty}$ for which

\begin{equation}
\ell(h) := \lim_{q\rightarrow\infty}\frac{1}{N_q}\sum_{n = 1}^{N_q}f(T^{n+h}x)\overline{y_n}\text{ exists for every }h \in \mathbb{N},
\end{equation}
we have

\begin{equation}
\lim_{H\rightarrow\infty}\frac{1}{H}\sum_{h = 1}^H|\ell(h^2)| = 0.
\end{equation}
\label{PolynomialProp}
\end{proposition}


We also see that Lemma \ref{keylemma} suggests we do not need the m.p.s. $([0,1],\mathscr{B},\mu,T)$ in Theorem \ref{ActualyFirstMainTheorem} to be weakly mixing. In particular, if we work with $L^2([0,1],\mu)$ instead of $L^1([0,1],\mu)$, then it seems that we only need $f \in L^2([0,1],\mu)$ to satisfy

\begin{equation}
\label{WeakMixing}
\lim_{N\rightarrow\infty}\frac{1}{N}\sum_{n = 1}^N|\langle U^nf,g\rangle-\int_Xfd\mu\int_Xgd\mu| = 0.
\end{equation}

\noindent for every $g \in L^2([0,1],\mu)$. However, it is not obvious as to whether or not this is the case. In the proof of Theorem $2.5$, we approximated $f$ by continuous functions, and every element of $C([0,1]) \subseteq L^2([0,1],\mu)$ satisfying equation \ref{WeakMixing} is the same as the m.p.s. $([0,1],\mathscr{B},\mu,T)$ being weakly mixing. This leads us to Conjecture \ref{HasToBeTrue}.


\begin{conjecture} Let $(X,\mathscr{B},\mu,T)$ be a m.p.s. and let $f \in L^2(X,\mu)$ satisfy

\begin{equation}
\lim_{N\rightarrow\infty}\frac{1}{N}\sum_{n = 1}^N|\langle U^nf, g\rangle| = 0
\end{equation}

\noindent for every $g \in L^2(X,\mu)$. Then $\left(f(T^nx)\right)_{n = 1}^{\infty}$ is a weakly mixing sequence for a.e. $x \in X$.
\label{HasToBeTrue}
\end{conjecture}


\section{A Pointwise Theorem for Strongly Mixing Systems}

\hskip 4mm In this section we prove Theorem \ref{mainStrongMixingTheorem}, which is Theorem \ref{SecondMainTheorem} for standard measure preserving systems. The method of extending Theorem \ref{mainStrongMixingTheorem} to \ref{SecondMainTheorem} is the same as the method used in section $4$ to extend \ref{FirstMainTheorem} to \ref{ActualyFirstMainTheorem} so we omit it to save space.

\begin{lemma} Let $([0,1],\mathscr{B},\mu,T)$ be a strongly mixing m.p.s. For every generic point $x \in [0,1]$ and each $f \in C([0,1])$ with $\int_Xfd\mu = 0$, $\left(f(T^nx)\right)_{n = 1}^{\infty}$ is a strongly mixing sequence.
\label{StrongMixingLemma}
\end{lemma}


\noindent{\it Proof.} Let $f \in C([0,1])$ be arbitrary and let $(y_n)_{n = 1}^{\infty} \subseteq \mathbb{C}$ be bounded. Let $(N_q)_{q = 1}^{\infty} \subseteq \mathbb{N}$ be any sequence for which

\begin{equation}
\lim_{q\rightarrow\infty}\frac{1}{N_q}\sum_{n = 1}^{N_q}f(T^{n+h}x)\overline{y_n}
\end{equation}

\noindent exists for every $h \in \mathbb{N}$. By Lemma \ref{keylemma}, let $g \in L^2([0,1],\mu)$ be such that

\begin{equation}
\lim_{q\rightarrow\infty}\frac{1}{N_q}\sum_{n = 1}^{N_q}f(T^{n+h}x)\overline{y_n} = \langle U^hf,g\rangle.
\end{equation}

\noindent Letting $U:L^2([0,1],\mu)\rightarrow L^2([0,1],\mu)$ denote the unitary operator induced by $T$, we see that $U$ is strongly mixing. It follows that

\begin{equation}
\pushQED{\qed}
\lim_{h\rightarrow\infty}\lim_{q\rightarrow\infty}\frac{1}{N_q}\sum_{n = 1}^{N_q}f(T^{n+h}x)\overline{y_n}h = \lim_{h\rightarrow\infty}\langle U^hf,g\rangle = 0. \qedhere
\popQED
\end{equation}


\begin{theorem}\label{mainStrongMixingTheorem} Let $([0,1],\mathscr{B},\mu,T)$ be a strongly mixing m.p.s. and let $f \in L^1([0,1],\mu)$ satisfy $\int_{[0,1]}fd\mu = 0$. For a.e. $x \in X$, $\left(f(T^nx)\right)_{n = 1}^{\infty}$ is a strongly mixing sequence.
\end{theorem}


\noindent{\it Proof.} Let $\epsilon > 0$ be arbitrary, and let $g \in C([0,1])$ be such that $||f-g||_1 \le \epsilon$ and $\int_Xgd\mu = 0$. Let $X \subseteq [0,1]$ be a set of full measure for which the Birkhoff Pointwise Ergodic Theorem holds for $|f-g|$ along every $x \in X$. We now see that

\begin{equation}
\lim_{N\rightarrow\infty}\frac{1}{N}\sum_{n = 1}^N|f(T^nx)-g(T^nx)| = \int_X|f-g|d\mu = ||f-g||_1 < \epsilon.
\end{equation}

Now let $(y_n)_{n = 1}^{\infty} \subseteq \mathbb{C}$ be uniformly bounded in norm by $1$. Since a.e. $x \in X$ is a generic point we may use Lemma \ref{WeakMixingForC(X)} to further refine $X$ to another set of full measure $X'$, such that for every $x \in X'$, $\left(g(T^nx)\right)_{n = 1}^{\infty}$ is a strongly mixing sequence. We see that for any $x \in X'$, we have

\begin{multline}
\lim_{h\rightarrow\infty}\overline{\lim_{N\rightarrow\infty}}|\frac{1}{N}\sum_{n = 1}^Nf(T^{n+h}x)\overline{y_n}| \\ \le \lim_{h\rightarrow\infty}\overline{\lim_{N\rightarrow\infty}}|\frac{1}{N}\sum_{n = 1}^Ng(T^{n+h}x)\overline{y_n}| +\epsilon = \epsilon.
\label{Modification2}
\end{multline}

\noindent Since $\epsilon > 0$ was arbitrary, we see that $\left(f(T^nx)\right)_{n = 1}^{\infty}$ is a strongly mixing sequence. \qed


In fact, we see that Theorem \ref{mainStrongMixingTheorem} can be used to characterize strongly mixing systems.


\begin{theorem} Let $(X,\mathscr{B},\mu,T)$ be a m.p.s. If for every $f \in L^{\infty}(X,\mu)$ with $\int_Xfd\mu = 0$ there exist a set $A_f \subseteq X$ satisfying $\mu(A_f) = 1$ and for every $x \in A_f$ we have that $\left(f(T^nx)\right)_{n = 1}^{\infty}$ is a strongly mixing sequence, then $T$ is strongly mixing.
\label{ActualStrongMixingConverse}
\end{theorem}


\noindent{\it Proof.} Let $A, B \in \mathscr{B}$ both be arbitrary. Let $X' \subseteq X$ be a set of full measure, such that for any $x \in X'$ and any $h \in \mathbb{N}$ we have that $\left(\mathbbm{1}_B(T^nx)-\mu(B)\right)_{n = 1}^{\infty}$, $\left(\mathbbm{1}_A(T^nx)-\mu(A)\right)_{n = 1}^{\infty}$ and $\left(\mathbbm{1}_{T^{-h}(A)\cap B}(T^nx)-\mu(T^{-h}(A)\cap B)\right)_{n = 1}^{\infty}$ are strongly mixing sequences. We note that if $(x_n)_{n = 1}^{\infty} \subseteq \mathbb{C}$ is a strongly mixing sequence, then for any $(N_q)_{q = 1}^{\infty} \subseteq \mathbb{N}$ for which the limits in equation \eqref{MixingImpliesErgodic} exist, we have

\begin{equation}
\label{MixingImpliesErgodic}
\lim_{N_q\rightarrow\infty}\frac{1}{N_q}\sum_{n = 1}^Nx_n = \lim_{h\rightarrow\infty}\lim_{q\rightarrow\infty}\frac{1}{N_q}\sum_{n = 1}^{N_q}x_{n+h}\overline{1} = 0.
\end{equation}

\noindent In particular, for any $x \in X'$, we have that

\begin{equation}
\lim_{N\rightarrow\infty}\frac{1}{N}\sum_{n = 1}^N\mathbbm{1}_B(T^nx) = \mu(B)
\end{equation}

\noindent and

\begin{equation}
\lim_{N\rightarrow\infty}\frac{1}{N}\sum_{n = 1}^N\mathbbm{1}_{T^{-h}(A)\cap B}(T^nx) = \mu(T^{-h}(A)\cap B)
\end{equation}

\noindent for every $h \in \mathbb{N}$. We now see that for any $x \in X'$, we have

\begin{equation}
\lim_{h\rightarrow\infty}\mu(T^{-h}(A)\cap B) = \lim_{h\rightarrow\infty}\lim_{N\rightarrow\infty}\frac{1}{N}\sum_{n = 1}^N\mathbbm{1}_{T^{-h}(A)\cap B}(T^nx)
\end{equation}

\begin{equation}
= \lim_{h\rightarrow\infty}\lim_{N\rightarrow\infty}\frac{1}{N}\sum_{n = 1}^N\mathbbm{1}_{T^{-h}(A)}(T^nx)\mathbbm{1}_{B}(T^nx)
\end{equation}

\begin{equation}
= \lim_{h\rightarrow\infty}\lim_{N\rightarrow\infty}\frac{1}{N}\sum_{n = 1}^N\mathbbm{1}_{A}(T^{n+h}x)\mathbbm{1}_{B}(T^nx)
\end{equation}

\begin{multline}
= \lim_{h\rightarrow\infty}\lim_{N\rightarrow\infty}\frac{1}{N}\sum_{n = 1}^N(\mathbbm{1}_{A}(T^{n+h}x)-\mu(A))\mathbbm{1}_{B}(T^nx)\\ +\lim_{h\rightarrow\infty}\lim_{N\rightarrow\infty}\frac{1}{N}\sum_{n = 1}^N\mu(A)\mathbbm{1}_{B}(T^nx)
\end{multline}

\begin{equation}
\pushQED{\qed}
= 0+\lim_{N\rightarrow\infty}\frac{1}{N}\sum_{n = 1}^N\mu(A)\mathbbm{1}_{B}(T^nx) = \mu(A)\mu(B). \qedhere
\popQED
\end{equation}


\section{Concluding Remarks}

\hskip 4mm The definitions of weakly mixing sequences and strongly mixing sequences are quite cumbersome since we have to pass to a subsequence $(N_q)_{q = 1}^{\infty}$ for which all of the relevant limits to exist. To circumvent this, we would like to propose definition \ref{SupremelyStronglyMixing} to replace the current definition of strongly mixing sequences.


\begin{definition} Let $(x_n)_{n = 1}^{\infty}$ be a sequence of complex numbers for which

\begin{equation}
\overline{\lim_{N\rightarrow\infty}}\frac{1}{N}\sum_{n = 1}^N|x_n| < \infty.
\end{equation}

$(x_n)_{n = 1}^{\infty}$ is a \noindent{\bf supremely strongly mixing sequence} if for any other bounded sequence of complex numbers $(y_n)_{n = 1}^{\infty}$ we have

\begin{equation}
\lim_{h\rightarrow\infty}\overline{\lim_{N\rightarrow\infty}}|\frac{1}{N}\sum_{n = 1}^Nx_{n+h}\overline{y_n}| = 0.
\end{equation}
\label{SupremelyStronglyMixing}
\end{definition}


It is clear that any supremely strongly mixing sequence is a strongly mixing sequence, and it would be convenient if the two notions were equivalent, as supremely strongly mixing suquences are much simpler to work with. However, Lemma \ref{CounterExample} shows us that there does not exist a supremely strongly mixing sequence.


\begin{lemma}\label{CounterExample} If $(x_n)_{n = 1}^{\infty}$ is any sequence of complex numbers for which

\begin{equation}
\overline{\lim_{N\rightarrow\infty}}\frac{1}{N}\sum_{n = 1}^N|x_n| < \infty,
\end{equation}

\noindent then there exists a bounded sequence of complex numbers $(y_n)_{n = 1}^{\infty}$ for which

\begin{equation}
\overline{\lim_{N\rightarrow\infty}}|\frac{1}{N}\sum_{n = 1}^Nx_{n+h}\overline{y_n}| = \overline{\lim_{N\rightarrow\infty}}\frac{1}{N}\sum_{n = 1}^N|x_n|,
\end{equation}

\noindent for every $h \in \mathbb{N}\cup\{0\}$.
\end{lemma}


\noindent{\it Proof.} Let $(N_k)_{k = 1}^{\infty}$ be such that

\begin{equation}
\overline{\lim_{N\rightarrow\infty}}\frac{1}{N}\sum_{n = 1}^N|x_n| = \lim_{k\rightarrow\infty}\frac{1}{N_k}\sum_{n = 1}^{N_k}|x_n|.
\end{equation}

\noindent By passing to a subsequence, we may assume without loss of generality that

\begin{equation}
\frac{1}{N_{k+1}}\left|\sum_{n = 1}^{N_k}|x_n|\right| < \frac{1}{k}.
\end{equation}

Let $f:\mathbb{N}_0^2\rightarrow\mathbb{N}_0$ be any bijection. For $N_{f(m,h)} < n \le N_{f(m,h)+1}$, let $y_n = \text{sgn}(x_{n+h})$. Now let $h \in \mathbb{N}_0$ be arbitrary, and note that

\begin{equation}
\overline{\lim_{N\rightarrow\infty}}\frac{1}{N}\sum_{n = 1}^N|x_n|
\end{equation}

\begin{equation}
\ge \overline{\lim_{N\rightarrow\infty}}|\frac{1}{N}\sum_{n = 1}^Nx_{n+h}\overline{y_n}| \ge \overline{\lim_{m\rightarrow\infty}}|\frac{1}{N_{f(m,h)+1}}\sum_{n = 1}^{N_{f(m,h)+1}}x_{n+h}\overline{y_n}|
\end{equation}

\begin{equation}
= \overline{\lim_{m\rightarrow\infty}}\frac{1}{N_{f(m,h)+1}}\left|\sum_{n = 1}^{N_{f(m,h)}}x_{n+h}\overline{y_n}+\sum_{n = N_{f(m,h)}+1}^{N_{f(m,h)+1}}|x_{n+h}|\right|
\end{equation}

\begin{equation}\ge \overline{\lim_{m\rightarrow\infty}}\frac{1}{N_{f(m,h)+1}}\left(\sum_{n = N_{f(m,h)}+1}^{N_{f(m,h)+1}}|x_{n+h}|-|\sum_{n = 1}^{N_{f(m,h)}}x_{n+h}\overline{y_n}|\right)
\end{equation}

\begin{equation}
\pushQED{\qed}
\ge \overline{\lim_{m\rightarrow\infty}}\frac{1}{N_{f(m,h)+1}}\sum_{n = 1}^{N_{f(m,h)+1}}|x_{n+h}|-\frac{2}{f(m,h)} = \overline{\lim_{N\rightarrow\infty}}\frac{1}{N}\sum_{n = 1}^N|x_n|. \qedhere
\popQED
\end{equation}


We will now show that Theorem \ref{ActualyFirstMainTheorem} and Theorem \ref{mainStrongMixingTheorem} hold for arbitrary measure preserving systems. Theorem \ref{ExtendingToTheGeneralSituation} is the generalization of Theorem \ref{ActualyFirstMainTheorem} to an arbitrary measure preserving system, and the argument for Theorem \ref{mainStrongMixingTheorem} is identical. For the proof of Theorem \ref{ExtendingToTheGeneralSituation} we will be using vocabulary and notation from chapter 15 of \cite{Royden} that we will not discuss here.


\begin{theorem} Let $(X,\mathscr{B},\mu,T)$ be a weakly mixing m.p.s. and let $f \in L^1(X,\mu)$ satisfy $\int_Xfd\mu = 0$. For a.e. $x \in X$, $\left(f(T^nx)\right)_{n = 1}^{\infty}$ is a weakly mixing sequence.
\label{ExtendingToTheGeneralSituation}
\end{theorem}


\noindent{\it Proof.} Let $\mathscr{A}$ denote the $\sigma$-algebra generated by $\{\{f < q\}_{q \in \mathbb{Q}}\}$. We see that $\mathscr{A}$ is a countably generated $\sigma$-algebra with respect to which $f$ is measurable. By Carath\'eodory's Theorem (cf. Theorem 15.3.4 of \cite{Royden}) there exists an isomorphism $\Phi$ of $<\mathscr{A}, \mu>$ into $<\mathscr{L}/\mathscr{N},m>$, where $\mathscr{L}$ is the Lebesgue $\sigma$-algebra on $[0,1]$, $m$ is Lebesgue measure and $\mathscr{N} \subseteq \mathscr{L}$ is the $\sigma$-algebra of null sets. By Proposition $15.2.2$ of \cite{Royden}, there exists $\tilde{f} \in L^1([0,1],m)$ be such that for every $q \in \mathbb{Q}$ we have $\{\tilde{f} < q\} = \Phi(\{f < q\})$. By Proposition $15.6.19$ of \cite{Royden}, let $\phi:X\rightarrow[0,1]$ be a measurable transformation for which $\mu(\phi^{-1}(B)\triangle\Phi^{-1}(B)) = 0$ for every $B \in \Phi(\mathscr{A})$. We see that for every $q \in \mathbb{Q}$, we have that $(\tilde{f}\circ\phi)^{-1}((-\infty,q)) = f^{-1}((-\infty,q))$, so $\tilde{f}\circ\phi = f$ a.e. by the uniqueness portion of Proposition $15.6.19$ of \cite{Royden}. Now note that $\Phi\circ T^{-1}\circ\Phi^{-1}$ is a $\sigma$-isomorphism of $\Phi(\mathscr{A})$ to itself, so yet another application of Proposition $15.6.19$ of \cite{Royden} yields a map $S:[0,1]\rightarrow[0,1]$ for which $S^{-1}(B) = \Phi(T^{-1}(\Phi^{-1}(B)))$ for every $B \in \Phi(\mathscr{A})$. Noting that $\phi^{-1}\circ S^{-1}$ and $T^{-1}\circ \phi^{-1}$ are the same $\sigma$-homomorphism, we can use the uniqueness portion of Proposition $15.6.9$ of \cite{Royden} once again to see that $S\circ \phi = \phi\circ T$. It follows that for some $X' \subseteq X$ with $\mu(X') = 1$ we have that $S^n(\phi(x)) = \phi(T^n(x))$ and $\tilde{f}(\phi(T^n(x))) = f(T^n(x))$ for all $n \in \mathbb{N}$ and all $x \in X'$. By Theorem $1$, let $X'' \subseteq [0,1]$ be such that $\left(\tilde{f}(S^nx)\right)_{n = 1}^{\infty}$ is a weakly mixing sequence for every $x \in X''$. Finally, we see that for any $x \in X'\cap\phi^{-1}(X'')$ we have that $\left(f(T^nx)\right)_{n = 1}^{\infty}$ is a weakly mixing sequence. \qed

\vskip 5mm

\begin{remark}
Lastly, we note that Theorems similar to Theorem \ref{ActualyFirstMainTheorem} and Theorem \ref{mainStrongMixingTheorem} can be proven for some other levels of mixing (such as mild mixing but not K-mixing) using the same techniques. We see that the only difference in the proofs of Theorem \ref{ActualyFirstMainTheorem} and Theorem \ref{mainStrongMixingTheorem} is the difference between equations \eqref{Modification} and \eqref{Modification2}. In particular, the mode of convergence in \eqref{Modification} and \eqref{Modification2} is the only difference. By working with convergence along filters a Theorem that simultaneously generalizes Theorems \ref{ActualyFirstMainTheorem} and \ref{mainStrongMixingTheorem} can be proven, and the author plans to include the statement and proof of this Theorem in his thesis. It is omitted here since it requires preliminary knowledge about convergence along filters.
\end{remark}

\subsection*{Acknowledgements}
I am grateful to Vitaly Bergelson for his guidance, support, and constructive criticism throughout the many drafts of this paper. I would also like to thank the referee for his helpful comments after reviewing the paper.

\end{document}